\documentclass[11pt]{article}
\usepackage[english]{babel}
\usepackage[latin1]{inputenc}
\usepackage{exscale}
\usepackage[centertags]{amsmath}
\usepackage{amsfonts,amstext,amssymb,amsthm}
\usepackage{comment}
\usepackage{newlfont}
\usepackage{enumerate}
\usepackage{graphicx}
\usepackage{these}
\usepackage{color}
\usepackage[usenames,dvipsnames]{xcolor}
\definecolor{BlueFonse}{rgb}{0,0,1}
\definecolor{BlueFonse1}{cmyk}{1,0,0,0.7}
\usepackage[colorlinks=true,linkcolor=blue,citecolor=ForestGreen]{hyperref}
\usepackage{url}
\title{The isoperimetric problem of a complete Riemannian manifold with a finite number of $C^0$-asymptotically Schwarzschild ends}
\author{Abraham Mu\~noz Flores\footnote{Partially supported by Capes}, Stefano Nardulli}
\begin{document}
      \maketitle
      \begin{center}
\noindent {\sc abstract}. We study the problem of existence of isoperimetric regions for large volumes, in $C^0$-locally asymptotically Euclidean Riemannian manifolds with a finite number of $C^0$-asymptotically Schwarzschild ends. Then we give a geometric characterization of these isoperimetric regions, extending previous results contained in \cite{EichmairMetzgerInv}, \cite{EichmairMetzgerJDG}, and \cite{BrendleEichmair}. Moreover strengthening a little bit the speed of convergence to the Schwarzschild metric we obtain existence of isoperimetric regions for all volumes for a class of manifolds that we named $C^0$-strongly asymptotic Schwarzschild, extending results of \cite{BrendleEichmair}. Such results are of interest in the field of mathematical general relativity. 
\bigskip\bigskip

\bigskip

\centerline{\bf AMS subject classification: }
49Q20, 58E99, 53A10, 49Q05.
\end{center}
      \tableofcontents       
      \newpage
\section{Introduction}\label{1}
 
The isoperimetric problem as defined in Definition \ref{Def:IsPWeak} is studied since the ancient times, its solution in the Euclidean plane and in the Euclidean $3$-dimensional space was known being the ball. Nevertheless only at the end of the nineteenth century, the first rigorous proof of this fact appeared. Since that time a lot of progress were made in the direction of proving existence and characterization of isoperimetric regions, but the list of manifolds for which we know the solution of the isoperimetric problem still remain too short. This paper is intended to expand this list. In general to prove existence and characterize geometrically isoperimetric regions is a quite hard task. In a series of papers \cite{EichmairMetzgerInv} (improving in various ways previous results in \cite{EichmairMetzgerJDG}), and \cite{BrendleEichmair}, M. Eichmair, J. Metzger, and S. Brendle consider the isoperimetric problem in a boundaryless initial data set $M$ that is also $C^0$-asymptotically  Scharzschild of mass $m>0$ with just one end. In this paper we consider the isoperimetric problem in a manifold with the same asymptotical conditions on the geometry but with a finite number of ends. Our characterization coincides with that of \cite{EichmairMetzgerInv} if the manifold have just one end, and with Corollary $16$ of \cite{BrendleEichmair} if $M$ is a double Schwarzschild manifold. To do this we use the theory developed in \cite{NarAsian}, \cite{NarCalcVar}, by both the two authors of this paper in \cite{FloresNardulli015}, the techniques developed in \cite{EichmairMetzgerInv}, \cite{EichmairMetzgerJDG}, and \cite{BrendleEichmair}. In particular we will use Theorem $4.1$ of \cite{EichmairMetzgerInv} and Theorem $15$ of \cite{BrendleEichmair}, combined with ad hoc new nontrivial arguments, needed to deal with the wider class of multiended manifolds $M$ considered here. The difficulties encountered to achieve the proof of the theorems are technical and they will become apparent later in the proofs.
\subsection{Finite perimeter sets in Riemannian manifolds}
We always assume that all the Riemannian manifolds $M$ considered are smooth with smooth Riemannian metric $g$. We denote by $V_g$ the canonical Riemannian measure induced on $M$ by $g$, and by $A_g$ the $(n-1)$-Hausdorff measure associated to the canonical Riemannian length space metric $d$ of $M$. When it is already clear from the context, explicit mention of the metric $g$ will be suppressed.  
\begin{Def} Let $M$ be a Riemannian manifold of dimension $n$, $U\subseteq M$ an open subset, $\mathfrak{X}_c(U)$ the set of smooth vector fields with compact support on $U$. Given $\Omega\subset M$ measurable with respect to the Riemannian measure, the \textbf{perimeter of $\Omega$ in $U$}, $ \mathcal{P}(\Omega, U)\in [0,+\infty]$, is
      \begin{equation}
                 \mathcal{P}(\Omega, U):=sup\left\{\int_{U}\chi_E div_g(X)dV_g: X\in\mathfrak{X}_c(U), ||X||_{\infty}\leq 1\right\},
      \end{equation}  
where $||X||_{\infty}:=\sup\left\{|X_p|_{{g}_p}: p\in M\right\}$ and $|X_p|_{{g}_p}$ is the norm of the vector $X_p$ in the metric $g_p$ on $T_pM$. If $\mathcal{P}(\Omega, U)<+\infty$ for every open set $U$, we call $E$ a \textbf{locally finite perimeter set}. Let us set $\mathcal{P}(\Omega):=\mathcal{P}(\Omega, M)$. Finally, if $\mathcal{P}(\Omega)<+\infty$ we say that \textbf{$\Omega$ is a set of finite perimeter}.    
\end{Def}
\begin{Def} We say that a sequence of finite perimeter sets $\Omega_j$ \textbf{converges in $L^1_{loc}(M)$} or in the \textbf{locally flat norm topology}, to another finite perimeter set $\Omega$, and we denote this by writing $\Omega_j\rightarrow\Omega$ in $L^1_{loc}(M)$, if $\chi_{\Omega_j}\rightarrow\chi_{\Omega}$ in $L^1_{loc}(M)$, i.e., if $V((\Omega_j\Delta \Omega)\cap U)\rightarrow 0,\;\forall U\subset\subset M$. Here $\chi_{\Omega}$ means the characteristic function of the set $E$ and the notation $U\subset\subset M$ means that $U\subseteq M$ is open and $\overline{U}$ (the topological closure of $U$) is compact in $M$.
\end{Def}   
\begin{Def}
We say that a sequence of finite perimeter sets $\Omega_j$ \textbf{converge in the sense of finite perimeter sets} to another finite perimeter set $\Omega$, if $\Omega_j\rightarrow\Omega$ in $L^1_{loc}(M)$, and 
\begin{eqnarray*} 
                           \lim_{j\rightarrow+\infty}\mathcal{P}(\Omega_j)=\mathcal{P}(\Omega).
\end{eqnarray*}
\end{Def}
For a more detailed discussion on locally finite perimeter sets and functions of bounded variation on a Riemannian manifold, one can consult \cite{MPPP}.
\subsection{Isoperimetric profile, compactness and existence of isoperimetric regions}
Standard results of the theory of sets of finite perimeter, guarantee that $A(\partial^*\Omega)=\mathcal{H}^{n-1}(\partial^*\Omega)=\mathcal{P}(\Omega)$ where $\partial^*\Omega$ is the reduced boundary of $\Omega$. In particular, if $\Omega$ has smooth boundary, then $\partial^*\Omega=\partial\Omega$, where $\partial\Omega$ is the topological boundary of $\Omega$. Furthermore, one can always choose a representative of $\Omega$ such that $\overline{\partial^*\Omega}=\partial\Omega$. In the sequel we will not distinguish between the topological boundary and the reduced boundary when no confusion can arise.
\begin{Def}\label{Def:IsPWeak}
Let $M$ be a Riemannian manifold of dimension $n$ (possibly with infinite volume). We denote by $\tilde{\tau}_M$ the set of  finite perimeter subsets of $M$. The function $\tilde{I}_M:[0,V(M)[\rightarrow [0,+\infty [$  defined by 
     $$\tilde{I}_M(v):= \inf\{\mathcal{P}(\Omega)=A(\partial \Omega): \Omega\in \tilde{\tau}_M, V(\Omega )=v \}$$ 
is called the \textbf{isoperimetric profile function} (or shortly the \textbf{isoperimetric profile}) of the manifold $M$. If there exists a finite perimeter set $\Omega\in\tilde{\tau}_M$ satisfying $V(\Omega)=v$, $\tilde{I}_M(V(\Omega))=A(\partial\Omega)= \mathcal{P}(\Omega)$ such an $\Omega$ will be called an \textbf{isoperimetric region}, and we say that $\tilde{I}_M(v)$ is \textbf{achieved}. 
\end{Def} 

Compactness arguments involving finite perimeter sets implies always existence of isoperimetric regions, but there are examples of noncompact manifolds without isoperimetric regions of some or every volumes. For further information about this point the reader could see the introduction of \cite{NarAsian} or \cite{MonNar} or Appendix $H$ of \cite{EichmairMetzgerInv} and the discussions therein. So we cannot have always a compactness theorem if we stay in a non-compact ambient manifold. If $M$ is compact, classical compactness arguments of geometric measure theory  combined with the direct method of the calculus of variations provide existence of isoperimetric regions in any dimension $n$. Hence, the problem of existence of isoperimetric regions in complete noncompact Riemannian manifolds is meaningful and in fact quite hard as we can argue from the fact that the list of manifolds for which we know whether isoperimetric regions exists or not, is very short. For completeness we remind the reader that if $n\leq 7$, then the boundary $\partial\Omega$ of an isoperimetric region is smooth. If $n\geq 8$ the support of the boundary of an isoperimetric region is the disjoint union of a regular part $R$ and a singular part $S$. $R$ is smooth at each of its points and has constant mean curvature, while $S$ has Hausdorff-codimension at least $7$ in $\partial\Omega$. For more details on regularity theory see \cite{Morg1} or \cite{Morgmt} Sect. $8.5$, Theorem $12.2$.  
 
\subsection{Main Results} 
The main result of this paper is the following theorem which is a nontrivial consequence of the theory developed in \cite{NarCalcVar}, \cite{NarAsian}, \cite{FloresNardulli}, \cite{FloresNardulli015}, combined with the work done in \cite{EichmairMetzgerInv}. This gives answers to some mathematical problems arising naturally in general relativity. 
\begin{Res}\label{Res:Existence} 
Let $(M^n,g)$ be an $n\geq 3$ dimensional complete boundaryless Riemannian manifold. Assume that there exists a compact set $K\subset\subset M$ such that $M\setminus K=\mathring{\bigcup}_{i\in\mathcal{I}}E_i$, where $\mathcal{I}:=\{1,...,l\}$, $l\in\N\setminus\{0\}$, and each $E_i$ is an end which is $C^0$-asymptotic to Schwarzschild of mass $m>0$ at rate $\gamma$, see Definition \ref{Def:AsymptoticSchwarzschild}. Then there exists $V_0=V_0(M,g)>0$ such that for every $v\geq V_0$ there exists at least one isoperimetric region $\Omega_v$ enclosing volume $v$. Moreover 
%$\partial\Omega_v$ 
%as exactly one connected component in exactly one end  smooth connected normal graph over the boundary $\partial B_r$, where $B_r$ is a coordinate ball with $V(B_r)=v$, i.e., 
$\Omega_v$ satisfies the conclusions of Lemma \ref{Lemma:IsopregionsLargeVolume}.
\end{Res}
\begin{figure}[h!]
  \centering
    \includegraphics[width=.9\linewidth]{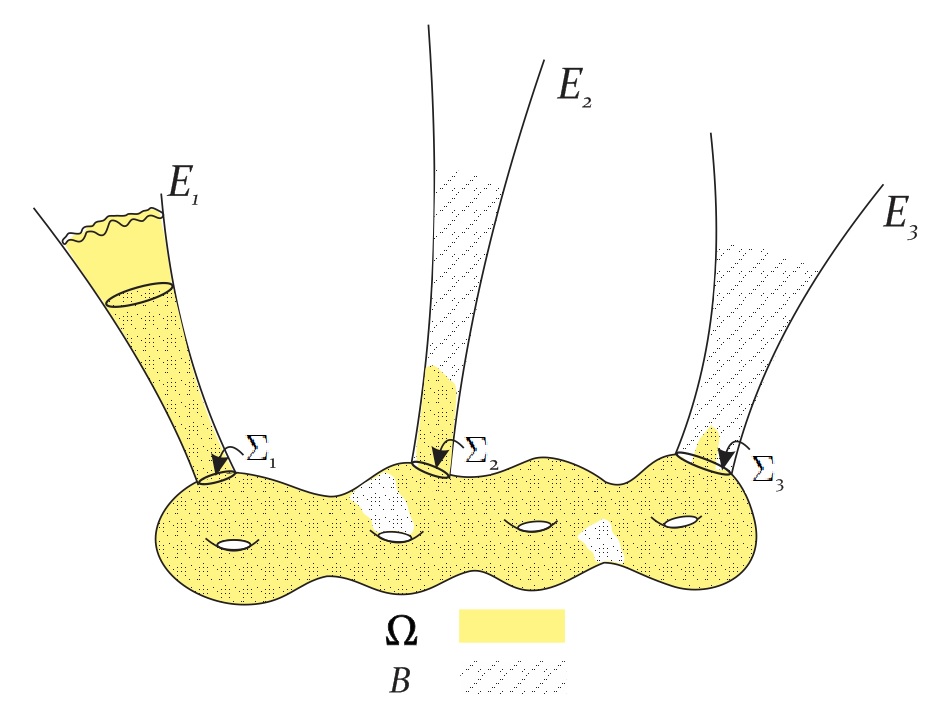}
  \caption{The isoperimetric region $\Omega$ is in yellow and $B$ is dotted.}
  \label{Fig:CharIsopRegLargeVolumes}
\end{figure}
\begin{Rem} The characterization of isoperimetric regions in Theorem \ref{Res:Existence} is achieved using Theorem $4.1$ of \cite{EichmairMetzgerInv} applied to the part of an isoperimetric region that have a sufficiently big volume in an end, the details and suitable modifications of the proof are presented in Lemma \ref{Lemma:IsopregionsLargeVolume}.
\end{Rem}
\begin{Rem}  Since the ends are like in \cite{EichmairMetzgerInv} it follows trivially that, if it happens that an end $E$ is $C^2$-asymptotic to Schwarzschild, then the volume $V_0$ can be chosen in such a manner that there exists a unique smooth isoperimetric (relatively to $E$) foliations of $E\setminus B$. Moreover, if $E$ is  asymptotically even (see Definition $2.1$ of \cite{EichmairMetzgerInv}) then the centers of mass of $\partial\Omega_v$ converge to the center of mass of $E$, as $V$ goes to $+\infty$, compare section $5$ of \cite{EichmairMetzgerInv}. 
\end{Rem}
\begin{CorRes}\label{Cor:multimass} If we allow $M$ in the preceding theorem to have each end $E_i$, with mass $m_i>0$. Then there exists a volume $V_0=V_0(M,g)>0$ and a subset $\mathcal{I}\subseteq\{1,..., \tilde{N}\}$, defined as $\mathcal{I}:=\{i:m_i=\max\{m_1,\dots, m_{\tilde{N}}\}\}$ such that for every volume $v\in [V_0, +\infty[$ there exist an isoperimetric region $\Omega_v$ that satisfies the conclusion of Lemma \ref{Lemma:IsopregionsLargeVolume} in which the preferred end $E_{\Omega_v}\in\{E_i\}_{i\in\mathcal{I}}$.  In particular, if $m_i\neq m_j$ for all $i\neq j$, then $\mathcal{I}=\{i\}$ is reduced to a singleton and this means that there exists exactly one end $E_i$ in which the isoperimetric regions for large volumes prefer to stay with a large amount of volume.
\end{CorRes}
In the next theorem paying the price of strengthening the rate of convergence to the Scwarzschild metric inside each end, we can show existence of isoperimetric regions in every volumes. The proof uses the generalized existence theorem of \cite{NarAsian} and a slight modification of the fine estimates for the area of balls that goes to infinity of Proposition $12$ of \cite{BrendleEichmair}. 
\begin{Res}\label{Res:Existence1} Let $(M^n,g)$ be an $n\geq 3$ dimensional complete boundaryless Riemannian manifold. Assume that there exists a compact set $K\subset\subset M$ such that $M\setminus K=\mathring{\bigcup}_{i\in\mathcal{I}}E_i$, where $\mathcal{I}:=\{1,...,l\}$, $l\in\N\setminus\{0\}$, and each $E_i$ is a $C^0$-strongly asymptotic to Schwarzschild of mass $m>0$ end, see Definition \ref{Def:AsymptoticSchwarzschild1}. Then for every volume $0<v<V(M)$ there exists at least one isoperimetric region $\Omega_v$ enclosing volume $v$. 
\end{Res}

%\subsection{Plan of the article}
%\begin{enumerate}
 %          \item  Section \ref{1} constitutes the introduction of the paper. We state the main results of the paper. 
 %          \item In Section \ref{3} we prove Theorems \ref{Res:Existence} and \ref{Res:Existence1}. 
%\end{enumerate}
\subsection{Acknowledgements}  
The authors would like to aknowledge Pierre Pansu, Andrea Mondino, Michael Deutsch, Frank Morgan for their useful comments and remarks. The first author wishes to thank the CAPES for financial support.         
\section{Proof of Theorems \ref{Res:Existence} and \ref{Res:Existence1}}\label{3} 
\subsection{Definitions and notations}
Let us start by recalling the basic definitions from the theory of convergence of manifolds, as exposed in \cite{Pet}. This will help us to state the main results in a precise way. 
\begin{Def} For any $m\in\mathbb{N}$, $\alpha\in [0, 1]$, a sequence of pointed smooth complete Riemannian manifolds is said to \textbf{converge in
the pointed $C^{m,\alpha}$, respectively $C^{m}$ topology to a smooth manifold $M$} (denoted $(M_i, p_i, g_i)\rightarrow (M,p,g)$), if for every $R > 0$ we can find a domain $\Omega_R$ with $B(p,R)\subseteq\Omega_R\subseteq M$, a natural number $\nu_R\in\mathbb{N}$, and $C^{m+1}$ embeddings $F_{i,R}:\Omega_R\rightarrow M_i$, for large $i\geq\nu_R$ such that $B(p_i,R)\subseteq F_{i,R} (\Omega_R)$ and $F_{i,R}^*(g_i)\rightarrow g$ on $\Omega_R$ in the $C^{m,\alpha}$, respectively $C^m$ topology. 
\end{Def}\noindent  
\begin{Def}\label{Def:BoundedGeometry}
A complete Riemannian manifold $(M, g)$, is said to have \textbf{bounded geometry} if there exists a constant $k\in\mathbb{R}$, such that $Ric_M\geq k(n-1)$ (i.e., $Ric_M\geq k(n-1)g$ in the sense of quadratic forms) and $V(B_{(M,g)}(p,1))\geq v_0>0$ for some positive constant $v_0$, where $B_{(M,g)}(p,r)$ is the geodesic ball (or equivalently the metric ball) of $M$ centered at $p$ and of radius $r> 0$.
\end{Def}
\begin{Rem} In general, a lower bound on $Ric_M$ and on the volume of unit balls does not ensure that the pointed limit metric spaces at infinity are still manifolds. 
\end{Rem} 
This motivates the following definition, that is suitable for most applications to general relativity for example. 
\begin{Def}\label{Def:BoundedGeometryInfinity}
We say that a smooth Riemannian manifold $(M^n, g)$ has $C^{k,\alpha}$-\textbf{locally asymptotic bounded geometry} if it is of bounded geometry and if for every diverging sequence of points $(p_j)$, there exist a subsequence $(p_{{j}_{l}})$ and a pointed smooth manifold $(M_{\infty}, g_{\infty}, p_{\infty})$ with $g_{\infty}$ of class $C^{k,\alpha}$ such that the sequence of pointed manifolds $(M, p_{{j}_{l}}, g)\rightarrow (M_{\infty}, g_{\infty}, p_{\infty})$, in  $C^{k,\alpha}$-topology.  
\end{Def}
For a more detailed discussion about this last definition and motivations, the reader could find useful to consult \cite{NarAsian}, we just make the following remark, illustrating some classes of manifolds for which Definition \ref{Def:BoundedGeometryInfinity} holds.
\begin{Rem} Observe that if $(M,g,p)\in\mathcal{M}^{k,\alpha}(n, Q, r)$ for every $p\in M$, then $M$ has $C^0$-bounded geometry. So Theorem 
\ref{Cor:Genexistence} applies to pointed manifolds in $\mathcal{M}^{k,\alpha}(n, Q, r)$. For the exact definitions see chapter 10 of \cite{Pet}. In particular according to M. Anderson (compare Theorem $76$ of \cite{Pet}) $M\in\mathcal{M}^{0,\alpha}(n, Q, r)$, whenever $|Ric_M|\leq\Lambda$, and $M$ has $V(B_{(M,g)}(p,1))\geq v_0>0$ for some positive constant $v_0$.
\end{Rem}  
\begin{Thm}[Generalized existence \cite{NarAsian}]\label{Cor:Genexistence}
Let $M$ have $C^0$-locally asymptotically bounded geometry. Given a positive volume $0<v < V(M)$, there are a finite number $N$, of limit manifolds at infinity such that their disjoint union with M contains an isoperimetric region of volume $v$ and perimeter $I_M(v)$. Moreover, the number of limit manifolds is at worst linear in $v$.
Indeed $N\leq\left[\frac{v}{v^*}\right]+1=l(n,k, v_0, v)$, where $v^*$ is as in Lemma 3.2 of \cite{EB}.
\end{Thm}
Now we come back to the main interest of our theory, i.e., to extend arguments valid for compact manifolds to noncompact ones. To this aim let us introduce
the following definition suggested by Theorem \ref{Cor:Genexistence}.
\begin{Def} We call $D_{\infty}=\bigcup_{i} D_{\infty, i}$ a finite perimeter set in $\tilde{M}$ a \textbf{generalized set of finite perimeter of $M$} and an isoperimetric region of $\tilde{M}$ a \textbf{generalized isoperimetric region}, where $\tilde{M}:=\mathring{\bigcup}_i M_{\infty, i}$, $D_{\infty, i}\subseteq M_{\infty, i}$ is isoperimetric in $M_{\infty, i}$.
\end{Def}
\begin{Rem} We remark that $D_{\infty}$ is a finite perimeter set of volume $v$ in $\mathring{\bigcup}_i M_{\infty, i}$. \end{Rem}
\begin{Rem} If $D$ is a genuine isoperimetric region contained in $M$, then $D$ is also a generalized isoperimetric region with $N=1$ and $$(M_{\infty,1}, g_{\infty,1})=(M,g).$$ This does not prevent the existence of another generalized isoperimetric region of the same volume having more than one piece at infinity.
\end{Rem}
\begin{Def}\label{Def:CmaAsFlat} Let $k\in\mathbb{N}$ and $\alpha\in [0,1]$ be given.  We say that a complete Riemannian $n$-manifold $(M,g)$ is \textbf{$C^{k,\alpha}$-locally asymptotically flat} or equivalently \textbf{$C^{k,\alpha}$-locally asymptotically Euclidean} if it is $C^{k,\alpha}$-locally asymptotic bounded geometry and for every diverging sequence of points $(p_j)_{j\in \mathbb{N}}$ there exists a subsequence $(p_{j_l})_{l\in \mathbb{N}}$ such that the sequence of pointed manifolds $$(M,g,p_{j_l})\to (\R^n, \delta, 0_{\R^n}),$$ in the pointed $C^{k,\alpha}$-topology, where $\delta$ is the canonical Euclidean metric of $\R^n$.
\end{Def}
%\begin{Rem} Observe that a $C^{k,\alpha}$-locally asymptotically Euclidean manifold in the sense of Definition \ref{Def:CmaAsFlat} is of bounded geometry in the sense of Definition \ref{Def:CmaAsFlat}, provided $k\geq2$.
%\end{Rem}
\begin{Def}\label{Def:Initialdataset} An \textbf{initial data set} $(M, g)$ is a connected complete
(with no boundary) $n$-dimensional Riemannian manifold such that there exists a positive constant $C>0$, a
bounded open set $U\subset M$, a positive natural number $\tilde{N}$, such that $M\setminus U=\mathring{\cup}_{i=1}^{\tilde{N}}E_i$, and $E_i\cong_{x_i} \R^n\setminus B_1(0)$, in the
coordinates induced by $x_i = (x_i^1, . . . , x_i^n)$ satisfying
\begin{equation}
r|g_{ij}-\delta_{ij}|+r^2|\partial_kg_{ij}|+r^3|\partial_{kl}^2g_{ij}|\leq C,
\end{equation}
 for all $r\geq 2$, where $r:=|x|=\sqrt{\delta_{ij}x^ix^j}$, (Einstein convention). We will use also the notations $B_r:=\{x\in\R^n:|x|<r\}$, and $S_r:=\{x\in\R^n:|x|=r\}$, $\tilde{B}^i_r:=\Phi_i(B_r)\subseteq M$ and $\tilde{S}^i_r:=\Phi_i(S_r)\subseteq M$, where $\Phi_i:=x_i^{-1}$, and we call $\tilde{B}^i_r$ a \textbf{centered coordinate ball of radius $r$} and $\tilde{S}^i_r$ a \textbf{centered coordinate sphere of radius $r$}, respectively. We put $\Sigma_i:=\tilde{S}^i_1$. Each $E_i$ is called an \textbf{end}. In what follows we suppress the index $i$ when no confusion can arise and we will note simply by $E$ an end, by $x$ the coordinate chart of $E$, $\tilde{B}_r:=\{p\in E:\:|x(p)|<r\}$, $\Sigma:=\{p\in E:\:|x(p)|=1\}$. 
\end{Def} 
\begin{Rem} An initial data set in the sense of Definition \ref{Def:Initialdataset} is $C^2$-locally asymptotically Euclidean in the sense of Definition \ref{Def:CmaAsFlat}.
\end{Rem}
\begin{figure}[h!]
  \centering
    \includegraphics[width=.7\linewidth]{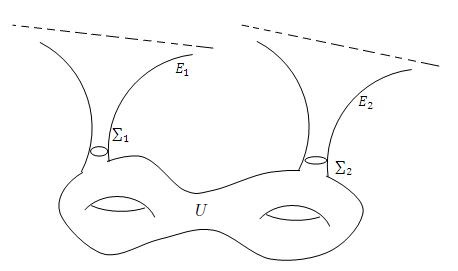}
  \caption{$M$ with $\Sigma_1=\partial E_1$ and $\Sigma_2=\partial E_2$, $E_1$ and $E_2$ ends.}
  \label{fig:fig2}
\end{figure}
\begin{figure}[h!]
  \centering
    \includegraphics[width=.5\linewidth]{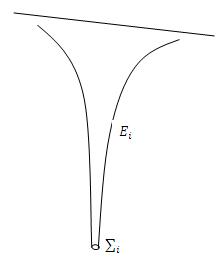}
  \caption{$\Sigma_i=\partial E_i$ and $E_{i}$ end.}
  \label{fig:fig2}
\end{figure}
In what follows we always assume that $n\geq 3$.
\begin{Def}\label{Def:AsymptoticSchwarzschild} For any $m>0$, $\gamma\in(0, 1]$, and $k\in\N$, we say that an initial data
set (compare Definition \ref{Def:Initialdataset}) is \textbf{$C^k$-asymptotic to Schwarzschild of mass $m>0$ at rate $\gamma$}, if
\begin{equation}\label{Eq:AsymptoticSchwarzschild1}
\sum_{l=0}^k r^{n-2+\gamma+l}|\partial^l(g-g_m)_{ij}|\leq C,
\end{equation}
for all $r\geq2$, in each coordinate chart $x_i:E_i\cong \R^n\setminus B_{\R^n}(0,1)$, where $(g_m)_{ij} = \left(1+ \frac{m}{2|x|^{n-2}}\right)^{\frac{4}{n-2}}\delta_{ij}$ is the usual Schwarzschild metric on $(\R^n\setminus\{0\})$. 
\end{Def}

\begin{figure}[h!]
\centering
 \includegraphics[width=.8\linewidth]{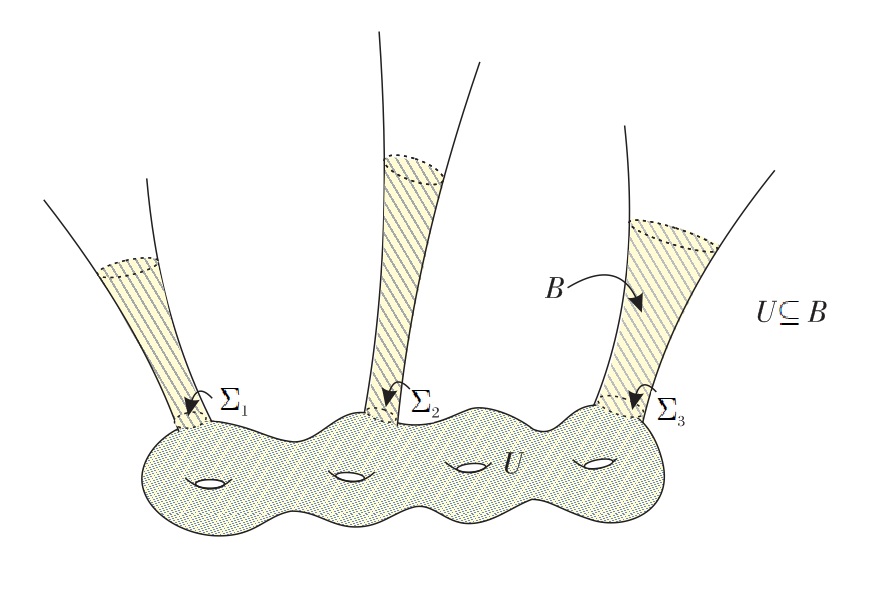}
  \caption{$M$ with $U\subseteq B$, $\Sigma_i=\partial E_i$ and $E_1$, $E_2$, $E_3$ ends.}
  \label{fig:fig5}
\end{figure}

\begin{Def}\label{Def:AsymptoticSchwarzschild1} For any $m>0$, $\gamma\in]0, +\infty[$, we say that an initial data
set is \textbf{$C^0$-strongly asymptotic to Schwarzschild of mass $m>0$ at rate $\gamma$}, if
\begin{equation}\label{Eq:AsymptoticSchwarzschild2}
 r^{2n+\gamma}|\left(g-g_m\right)_{ij}|\leq C,
\end{equation}
for all $r\geq2$, in each coordinate chart $x_i:E_i\cong \R^n\setminus B_{\R^n}(0,1)$, where $(g_m)_{ij} = \left(1+ \frac{m}{2|x|^{n-2}}\right)^{\frac{4}{n-2}}\delta_{ij}$ is the usual Schwarzschild metric on $(\R^n\setminus\{0\})$.
\end{Def}
\begin{Rem} If $M$ satisfies the assumptions of the Definition \ref{Def:AsymptoticSchwarzschild} or Definition \ref{Def:AsymptoticSchwarzschild1}, then $M$ is trivially $C^2$-locally asymptotically Euclidean in the sense of Definition \ref{Def:CmaAsFlat}, since $M$ is assumed to be an initial data set.
\end{Rem}
\begin{Lemme}\label{Lemma:IsopregionsLargeVolume}There exists $V_0=V_0(M,g)>0$, and a large ball $B$ such that if $\Omega\subseteq M$ is an isoperimetric region with $V(\Omega)=v\geq V_0$, then there exists an end $E_i=E_{\Omega}$ such that $\Omega\cap E_{\Omega}$ is the region below a normal graph based on $\tilde{S}^i_r$ where $V_g(\Omega\cap E_{\Omega})=V_g(\tilde{B}_r)$, i.e., $\Omega=x_i^{-1}(\varphi(B_r\setminus B_1))\mathring{\cup}\Omega^*$, with $\Omega^*\subseteq B$ and $\varphi(B_r\setminus B_1)\subseteq\R^n\setminus B_{\R^n}(0,1)$ is a suitable perturbation of $B_r\setminus B_1$. $\Omega\cap E_{\Omega}$ contains $\Sigma_i$ and is an isoperimetric region as in Theorem 4.1 of \cite{EichmairMetzgerInv}, $\Omega\setminus E_{\Omega}$ contains $\Sigma_i$ and $\Omega\setminus E_{\Omega}$ has least relative perimeter with respect to all domains in $B\setminus E$ containing $\Sigma_i$ and having volume equal to $V(\Omega\setminus E_{\Omega})$.
\end{Lemme}
\begin{Rem} In general $B$ contains $U$ and is much larger than $U$, see figure \ref{fig:fig5}. $B$ could be chosen in such a way that $M\setminus B$ is a union of ends that are foliated by the boundary of isoperimetric regions of that end, provided this foliation exists. Furthermore $B$ contain $U$ and all the $\tilde{B}^i_r$, with $r$ large enough to enclose a volume bigger than the volume $V_0$ given by Theorem $4.1$ of \cite{EichmairMetzgerInv}. 
\end{Rem}
\begin{figure}[h!]
  \centering
    \includegraphics[width=.7\linewidth]{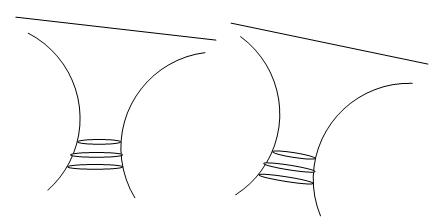}
  \caption{$M\setminus B$}
  \label{fig:fig2}
\end{figure}
\begin{Rem} Corollary $16$ of \cite{BrendleEichmair} is a particular instance of Lemma \ref{Lemma:IsopregionsLargeVolume} when the number of ends is two. Of course, in Corollary $16$ of \cite{BrendleEichmair} more accurate geometrical informations are given due to the very special features of the  double Schwarzschild manifolds considered there. See figure \ref{fig:fig4} in which the same notation of Corollary $16$ of \cite{BrendleEichmair} are used.
\end{Rem}
\begin{figure}[h!]
\centering
 \includegraphics[width=.6\linewidth]{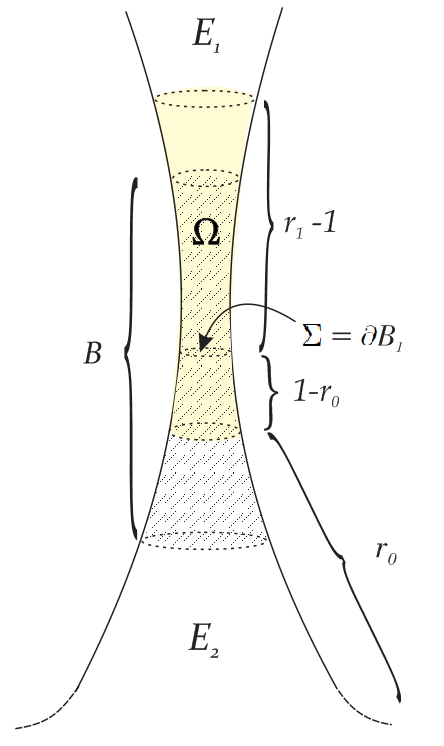}
  \caption{In Corollary $16$, of \cite{BrendleEichmair}, $(M,g)$ is the double Schwarzschild manifold of positive mass $m>0$, $(\R^n\setminus 0, g_m)$. The isoperimetric region $\Omega$ is in yellow, $B$ is dotted, $E_1$ is the preferred end. Just in this example $\Sigma$ is outermost minimal, w.r.t. any end.}
  \label{fig:fig4}
\end{figure}
\begin{Dem} In first observe that the existence of a geodesic ball $B$ satisfying the conclusions of the Lemma is essentially equivalent to the weaker assumption that the same geodesic ball $B$ contain just all the volume of $\Omega\setminus E$. To see this is trivial because we know that we can modify a finite perimeter set on a set of measure zero and stay always in the same equivalence class. Now, assume that a geodesic ball  $B$ satisfying the conclusions of Lemma \ref{Lemma:IsopregionsLargeVolume} does not exists. Let $(\Omega_j)$ be a sequence of isoperimetric regions in $M$, such that $V(\Omega_j)\rightarrow +\infty$. It follows easily using the fact that the number of ends is finite that there exists at least one end $E_{i_j}=:E_j$, such $V(\Omega_j\cap E_j)\rightarrow+\infty$. The crucial point is to show that  this end $E_{i_j}$ is unique. To show this we observe that the proof of Theorem $4.1$ of \cite{EichmairMetzgerInv}, applies exactly in the same way to our sequence $(\Omega_j)$ and our manifold $M$. This application of Theorem $4.1$ of \cite{EichmairMetzgerInv} gives us a volume $V_0$ that depends only on the geometric data of the ends such that if $\Omega$ is an isoperimetric region of volume $v\geq V_0$ then there exist an end $E$ such that $\Omega\cap E$ contains a large centered ball $\tilde{B}_{r/2}$ with $V(\tilde{B}_r)=V(\Omega\cap E)$. In particular, this discussion shows that for large values of the enclosed volume $v\geq V_0$ an isoperimetric region $\Omega$ is such that
\begin{equation}\label{Eq:IsopregionsLargeVolume}
\Sigma\subseteq \Omega,
\end{equation} 
and finally that $\Omega\cap E=x^{-1}(\varphi\left(B_r\setminus B_1\right))$,  
for larges values of $V(\Omega)$, where $\Sigma$ is the boundary of $E$. 
%Without loss of generality we can assume that at least for one end $E$ we have $V(\Omega_j\cap E)\rightarrow+\infty$, where $E$ is a fixed end with respect to $j$. 
Now we show that there is no infinite volume in more than one end. Roughly speaking, this follows quickly from the estimates (\ref{Eq:Main3-}) (a particular case of which is $(3)$ of \cite{EichmairMetzgerJDG}) because the dominant term in the expansion of the area with respect to volume is as the isoperimetric profile of the Euclidean space which shows that two big different coordinate balls each in one different end do worst than one coordinate ball of  in just one end of the the volume that the sum of the other two balls. To show this rigorously we will argue by contradiction. Assume that for every $j$, there exist two distinct ends $E_{j1}\neq E_{j2}$ such that $V(\Omega_j\cap E_{j1})=:v_{j1}\rightarrow+\infty$ and $V(\Omega_j\cap E_{j2})=:v_{j2}\rightarrow+\infty$. Again an application of Theorem 4.1 of  \cite{EichmairMetzgerInv} permits us to say that $\Omega_j\cap E_{j1}$ and $\Omega_j\cap E_{j2}$ are perturbations of large coordinates balls whose expansion of the area with respect to the enclosed volume is given by (\ref{Eq:Main3-}). Now put $\Omega'_j:=\Omega_j\setminus E_{j2}\cup\Omega''_j$, in such a way $\Omega'_j\cap E_{j1}$ is isoperimetric in $E_{j1}$ and $V(\Omega'_j)=V(\Omega_j)=v$, we have that
\begin{equation}
A(\partial\Omega'_j)=I_{\R^n}(v_{j1}+v_{j2})-m^*\left(v_{j1}+v_{j2}\right)^{\frac{1}{n}}+A(\Sigma_2)+\cdots.
\end{equation}
\begin{equation}
A(\partial\Omega_j)=I_{\R^n}(v_{j1})+I_{\R^n}(v_{j2})-m^*v_{j1}^{\frac{1}{n}}-m^*v_{j2}^{\frac{1}{n}}+\cdots.
\end{equation}
\begin{eqnarray}
A(\partial\Omega'_j)-A(\partial\Omega_j) & = & I_{\R^n}(v_{j1}+v_{j2})-I_{\R^n}(v_{j1})+I_{\R^n}(v_{j2})\\
& - & m^*\left(v_{j1}+v_{j2}\right)^{\frac{1}{n}}+m^*v_{j1}^{\frac{1}{n}}+m^*v_{j2}^{\frac{1}{n}}+\cdots\\
& \to & -\infty, 
\end{eqnarray}
when $v_{j1},v_{j2}\rightarrow+\infty$. Which contradicts the hypothesis that $\Omega_j$ is a sequence of isoperimetric regions.  We have to show that the diameters of $\Omega_j\cap(M\setminus E_j)$ are uniformly bounded w.r.t. $j$. 
We start this arguments by noticing that $V(\Omega_j\cap(M\setminus E_j))\leq K$ are uniformly bounded, for the same fixed positive constant $K$ and so we can pick a representative of $\Omega_j$ such that $\Omega_j\cap(M\setminus E_j)$ is entirely contained in $B$. Analogously to what is done in the proof of Theorem $3$ of \cite{NarAsian} or Lemma $3.8$ of \cite{Nar07} (these proofs were inspired by preceding works of Frank Morgan \cite{Morgan94} proving boundedness of isoperimetric regions in the Euclidean setting and Manuel Ritor\'e and Cesar Rosales in Euclidean cones \cite{RRosales}) we can easily conclude that 
\begin{equation} 
diam(\Omega_j\cap(M\setminus E_j))\leq C(n,k, v_0)V(\Omega_j\cap(M\setminus E_j))^{\frac{1}{n}}\leq CK^{\frac{1}{n}}, 
\end{equation}
which ensures the existence of our big geodesic ball $B$.
According to (\ref{Eq:IsopregionsLargeVolume}) we have $\Sigma\subset\Omega\setminus E_{\Omega}$, to finish the proof we prove that $\Omega\setminus E_{\Omega}\subseteq B\setminus E_{\Omega}$ is such that $A(\partial\Omega\cap (B\setminus E_{\Omega}))$ is equal to
\begin{equation}\label{Eq:IsopregionsLargeVolume1}
\inf\left\{A(\partial D\cap (B\setminus E_{\Omega})): D\subseteq B\setminus E_{\Omega}, \Sigma\subseteq D, V(D)=V(\Omega\setminus E_{\Omega})\right\}.
\end{equation}
Now we proceed to the detailed verification of (\ref{Eq:IsopregionsLargeVolume1}). In fact if (\ref{Eq:IsopregionsLargeVolume1}) was not true we can find a finite perimeter set $\Omega'$ (that can be chosen open and bounded with smooth boundary, but this does not matter here) inside $B\setminus E$ such that $\Sigma\subseteq\Omega'$, $V(\Omega')=V(\Omega\setminus E_{\Omega})$, and $$A(\partial\Omega')<A(\partial\Omega\cap (B\setminus E_{\Omega})),$$ but if this is the case we argue that $V(\left(\Omega\cap E_{\Omega}\right)\cup\Omega')=v$ and $$A\left(\partial\left[(\Omega\cap E_{\Omega}\right)\cup\Omega'\right])<A(\partial\Omega),$$ which contradicts the fact that $\Omega$ is an isoperimetric region of volume $v$.  
\end{Dem}

Now we prove Theorem \ref{Res:Existence}. 
%In this proof we use in a crucial way the existence of $B$ like in Lemma \ref{Lemma:IsopregionsLargeVolume}.

\begin{Dem}
      Take a sequence of volumes $v_i\rightarrow+\infty$. Applying the generalized existence Theorem $1$ of \cite{NarAsian}, we get that there exists
      $\Omega_i\subset M$, ($\Omega_i$ is eventually empty) isoperimetric region with $V(\Omega_i)=v_{i1}$ and $B_{\R^n}(0,r_i)\subset\R^n$ with $V(B_{\R^n}(0,r_i))=v_{i2}$, satisfying $v_{i1}+v_{i2}=v_i$, and $I_M(v_i)=I_M(v_{i1})+I_M(v_{i2})$. We observe that $I_M(v_{i1})=A(\partial\Omega_i)$ and that we have just one piece at infinity because two balls do worst than one in Euclidean space. Note that this argument was already used in the proof of Theorem $1$ of \cite{MonNar}. If $v_{i2}=0$ there is nothing to prove, the existence of isoperimetric regions follows immediately. If $v_{i2}>0$ one can have three cases
      \begin{enumerate}
                \item $v_{i1}\rightarrow+\infty$, 
                \item there exist a constant $K>0$ such that $0<v_{i1}\leq K$ for every $i\in\N$,
                \item  $v_{i1}=0$, for $i$ large enough.
      \end{enumerate}          
                We will show, in first, that we can rule out case 2) and 3). To do this, suppose by contradiction that $0<v_{i1}\leq K<+\infty$ then remember that by Theorem 1 of 
      \cite{FloresNardulli} the isoperimetric profile function $I_M$ is continuous so $V(\Omega_i)+A(\partial\Omega_i)\leq K_1$ where $K_1>0$ is another positive constant. We can extract from the sequence of volumes $v_{i1}$ a convergent subsequence named again $v_{i1}\rightarrow\bar{v}\geq 0$. By generalized existence we obtain a generalized isoperimetric region $D\subset \tilde{M}$ such that $V(D)=\bar{v}$, $I_M(\bar{v})=A(\partial D)$. Again $D=D_1\mathring{\cup}D_{\infty}$, with $D_1\subset M$ and $D_{\infty}\subset\R^n$ isoperimetric regions in their respective volumes and in their respective ambient manifolds. Hence $D_{\infty}$ is an Euclidean ball. But also by the continuity of $I_M=I_{\tilde{M}}$ we have that $D\mathring{\cup}B_{\R^n}(0,r_i)$ is a generalized isoperimetric region of volume $\bar{v}+v_{i2}$, it follows that $H_{\partial D_{\infty}}=\frac{n-1}{r_i}$ for every $i\in\N$. As a consequence of the fact that $v_i\rightarrow+\infty$ and $(v_{i1})$ is a bounded sequence we must have $v_{i2}\rightarrow+\infty$, hence $r_i\rightarrow+\infty$, and we get $H_{\partial D_{\infty}}=\lim_{r_i\rightarrow+\infty}\frac{n-1}{r_i}=0$. As it is easy to see it is impossible to have an Euclidean ball with finite positive volume and zero mean curvature. This implies that $D_{\infty}=\emptyset$, for $v_i$ large enough. As a consequence of the proof of Theorem $2.1$ of \cite{RRosales} or Theorem $1$ of \cite{NarAsian} and the last fact we have $\Omega_i\rightarrow D_1$ in the sense of finite perimeter sets of $M$. This last assertion implies that $V(D_1)=\bar{v}=\lim_{i\rightarrow+\infty} v_{i1}$. By Lemma 2.7 of \cite{NarAsian} we get $I_M\leq I_{\R^n}$. It follows that 
\begin{equation}\label{Eq:Main0}
          0\leq I_M(v_{i1})\leq I_{\R^n}(v_i)-I_{\R^n}(v_{i2})\rightarrow 0,
\end{equation}          
      because $v_i-v_{i2}\rightarrow\bar{v}$ and $I_{\R^n}$ is the function $v\mapsto v^{\frac{n-1}{n}}$, with fractional exponent $0<\frac{n-1}{n}<1$. By (\ref{Eq:Main0}) $\lim_{i\rightarrow+\infty} I_M(v_{i1})=0$, since $I_M$ is continuous we obtain $$\lim_{i\rightarrow+\infty} I_M(v_{i1})=I_M(\bar{v})=0=A(\partial D_1),$$ which implies that $V(D_1)=\bar{v}=0$. Now for small nonzero volumes, isoperimetric regions are psedobubbles with small diameter and big mean curvature $H_{\partial\Omega_i}\rightarrow +\infty$, because $M$ is $C^2$-locally asymptotically Euclidean, compare \cite{NarCalcVar} (for earlier results in the compact case compare \cite{NarAnn}),  but this is a contradiction because by first variation of area $H_{\partial\Omega_i}=H_{\partial B_{\R^n}(0,r_i)}=\frac{n-1}{r_i}$, with $r_i\rightarrow+\infty$.  We have just showed that $v_{i1}=0$ for $i$ large enough provided $v_{i1}$ is bounded, that is case 2) is simply impossible.
                 
Consider, now, the case 3), i.e., $v_{i1}=0$. To rule out this case we compare a large Euclidean ball of enclosed volume $v_{i2}$ with $\Omega_v:=x_i^{-1}(B_r\setminus B_1)$ choosing $r$ such that $V(\Omega_{v_{i2}})=v_{i2}$, by $(\ref{Eq:Main2})$, we get $A(\partial\Omega_{v_{i2}})\leq c_nv_{i2}^{\frac{n-1}{n}}$. If $v_{i1}=0$, for large $i$ then we have that all the mass stays in a manifold at infinity and so if we want to have existence we need an isoperimetric comparison for large volumes between $I_M(v)$ and $I_{M_{\infty}}(v)=I_{\R^n}(v)$. This isoperimetric comparison is a consequence of (\ref{Eq:Main2}) which gives that there exists a volume $v_0=v_0(C, m)$ (where $C$ is as in Definition \ref{Def:AsymptoticSchwarzschild}) such that 
\begin{equation}\label{Eq:Main1}
I_M(v)<I_{\R^n}(v),
\end{equation} for every $v\geq v_{0}$.  
To see this we look for finite perimeter sets $\Omega'_v\subset M$ which are not necessarily isoperimetric regions, which have volume $V(\Omega'_v)=v$ and $A(\partial\Omega'_v)<I_{\R^n}(v)$. A candidate for this kind of domains are coordinate balls inside a end $\tilde{B}_r:=x_i^{-1}(B_r\setminus B_1(0))$, with $r$ such that $V(\tilde{B}_r)=v$, because after straightforward calculations 
\begin{equation}\label{Eq:Main2} 
A(\partial\tilde{B}_{r(v)})=I_{\R^n}(v)-m^*v^{\frac{1}{n}}+o(v^{\frac{1}{n}})=c_nv^{\frac{n-1}{n}}-m^*v^{\frac{1}{n}}+o(v^{\frac{1}{n}}),
\end{equation} 
where $m^*>0$ is the same coefficient that appears in the asymptotic expansion of $$A_{g_m}(\partial\Omega_{v_m})=I_{\R^n}(v_m)-m^*v_m^{\frac{1}{n}}+o(v^{\frac{1}{n}}),$$ $v_m:=V_{g_m}(\Omega_v)$. Namely $m^*=c'_nm>0$, where $c'_n$ is a dimensional constant that depends only on the dimension $n$ of $M$. The calculation of $m^*$ is straightforward and we omit here the details, in the case of $n=3$ it comes immediately from (3) of \cite{EichmairMetzgerJDG}. It is worth to note here that the assumption (\ref{Eq:AsymptoticSchwarzschild1}) in Definition \ref{Def:AsymptoticSchwarzschild}, is crucial to have the remainder in (\ref{Eq:Main2}) of order of infinity strictly less than $v^{\frac{1}{n}}$. If the rate of convergence of $g$ to $g_m$ was of the order $r^{-\alpha}$ with $0<\alpha\leq n-2$ then this could add some extra term to $m^*$ in the asymptotic expansion (\ref{Eq:Main2}) that we could not control necessarily. This discussion permits to exclude case 3).

So we are reduced just to the case 1).
We will show that the only possible phenomenon that can happen is $v_{i1}\rightarrow +\infty$ and  $v_{i2}=0$. With this aim in mind we will show that it is not possible to have $v_{i1}\rightarrow+\infty$ and also $v_{i2}\rightarrow+\infty$ at the same time. A way to see this fact is to consider equation (\ref{Eq:Main2}) and observe that the leading term is Euclidean, now we take all the mass $v_{i2}$ and from infinity we add a volume $v_{i2}$ to the part in the end $E$, in this way we construct a competitor set (as in the proof of Lemma \ref{Lemma:IsopregionsLargeVolume}) $\tilde{\Omega}_{v_i}$ which is isoperimetric in the preferred end $E_i$ and such that $\tilde{\Omega}_{v_i}\setminus\Omega_{v_{i1}}=x_i^{-1}(\varphi(B_{r_i}\setminus B_1))$, where $E$ is one fixed end in which $V(E\cap\Omega_{v_{i1}})\rightarrow +\infty$, $r_i>1$, and $V(x_i^{-1}(\varphi(B_{r_i}\setminus B_1)))=v_{i2}$ in such a way that $\tilde{\Omega}_{v_i}\cap E$ is an isoperimetric region containing $\Sigma$ of $E$, i.e., a pertubation of a large coordinate ball as prescribed by Theorem 4.1 of \cite{EichmairMetzgerInv}, $\tilde{\Omega}_{v_i}\cap (B\setminus E)=:\bar{\Omega}_{\bar{v}_i}$ and $V(\tilde{\Omega}_{v_i})=v_i$. Hence by virtue of (\ref{Eq:Main2})  we get for large $v_{i1}$
\begin{equation}\label{Eq:Main3-}
I_M(v_{i1})=A(\partial B_{r(v_{i1})})+c(\bar{v}_{i1})+\varepsilon(v_{i1})=I_{\R^n}(v_{i1})-m^*v_{i1}^{\frac{1}{n}}+o(v_{i1}^{\frac{1}{n}}),
\end{equation}
where $\varepsilon(v_{i1})\rightarrow 0$ when $v_{i1}\rightarrow +\infty$, $c(\bar{v}_{i1})$ is the relative area of the isoperimetric region $\bar{\Omega}_{\bar{v}_{i1}}$ of volume $\bar{v}_{i1}$ inside $B\setminus E$ where $B$ is the fixed big ball of Lemma \ref{Lemma:IsopregionsLargeVolume}. It is easy to see that $\Omega_{v_{i1}}\cap \left(B\setminus E\right)$ could be caracterized as the isoperimetric region for the relative isoperimetric problem in $B\setminus E$ which contain the boundary $\Sigma$ of $E$. Such a relative isoperimetric region $\Omega'$ exists by standard compactness arguments of geometric measure theory, and regularity theory as in \cite{EichmairMetzgerJDG}, (compare also Theorem $1.5$ of \cite{DuzaarSteffen92}), in particular $A(\partial\Omega' \cap (B\setminus E))$ is equal to
\begin{equation}
\inf\left\{A(\partial D\cap (B\setminus E)): D\subseteq B, \Sigma\subseteq D, V(D)=V(\Omega'\setminus E)\right\}.
\end{equation}
Again by compactness arguments it is easy to show that the relative isoperimetric profile $I_{B\setminus E}:[0, V(B\setminus E)]\rightarrow [0, +\infty[$ is continuous (one can see this using the proof Theorem $1$ of \cite{FloresNardulli} that applies because we are in bounded geometry), and so $||I_{B\setminus E}||_{\infty}=c<+\infty$. If one prefer could rephrase this in terms of a relative Cheeger constant. This shows that $c(v)\leq c$ for every $v\in[0, V(B\setminus E)]$. This last fact legitimate the second equality in equation (\ref{Eq:Main3-}).
Thus readily follows
\begin{equation}\label{Eq:Main3}  
A(\partial\tilde{\Omega}_{v_i})<I_M(v_{i1})+I_{\R^n}(v_{i2})=I_M(v_i),
\end{equation}
 for large volumes $v_{i}\rightarrow+\infty$, which is the desired contradiction. We remark that the use Lemma \ref{Lemma:IsopregionsLargeVolume} is crucial to have the right shape of $\Omega_{v_{i1}}$ inside the preferred end $E$. To finish the proof, the only case that remains to rule out is when $v_{i1}\rightarrow +\infty$ and $0<v_{i2}\leq const.$ for every $i$. By the generalized compactness Theorem $1$ of \cite{FloresNardulli015} there exists $v_2\geq 0$ such that $v_{i2}\rightarrow v_2$. If $v_2>0$ then comparing the mean curvatures like already did in this proof, to avoid case 2) we obtain a contradiction, because the mean curvature of a large coordinate sphere tends to zero but the curvature of an Euclidean ball of positive volume $v_2$ is not zero. A simpler way to see this is again to look at formula (\ref{Eq:Main2}), since the leading term is $I_{\R^n}$ that is strictly subadditive, we can consider again a competing domain $\tilde{\Omega}_{v_i}$ such that $\tilde{\Omega}_{v_i}\setminus\Omega_{v_{i1}}=x^{-1}(\varphi(B_{r_i}\setminus B_1))$, with $E_i$ is such that $V(E_{i}\cap\Omega_{v_{i1}})\rightarrow +\infty$, $r_i>1$, and $V(x^{-1}(\varphi(B_{r_i}\setminus B_1)))=v_{i2}$, (\ref{Eq:Main3}) implies the claim.     If $v_2=0$ the situation is even worst because the mean curvature of Euclidean balls of volumes going to zero goes to $+\infty$, again because isoperimetric regions for small volumes are nearly round ball, i.e., pseudobubbles as showed in \cite{NarCalcVar}, whose theorems apply here since $M$ is $C^2$-locally asymptotically Euclidean. Hence we have necessarily that for $v_i$ large enough $v_{i2}=0$, which implies existence of isoperimetric regions of volumes $v_i$, provided $v_i$ is large enough. Since the sequence $v_i$ is arbitrary the first part of the theorem is proved. Now that we have established existence of isoperimetric regions for large volumes.      
   The second claim in the statement of Theorem \ref{Res:Existence} follows readily from Lemma \ref{Lemma:IsopregionsLargeVolume}. 
   %redo the same arguments of  Theorem 4.1 \cite{EichmairMetzgerInv}, \cite{EichmairMetzgerJDG} (which uses results from \cite{Ros}, \cite{Nar07}, \cite{MJ} combined with a highly nontrivial effective isoperimetric inequality, Lemma 3.4 of \cite{EichmairMetzgerInv}) applied to the part of the isoperimetric region that is inside the selected end within a big volume as did in Lemma \ref{Lemma:IsopregionsLargeVolume}.
      \end{Dem}    
\begin{Rem} If we allow to each end $E_i$ of $M$ to have a mass $m_i>0$ that possibly is different from the masses of the others ends, then we can guess in which end the isoperimetric regions for big volumes concentrates with "infinite volume". In fact the big volumes isoperimetric regions will prefer to stay in the end that for big volumes do better isoperimetrically and by (\ref{Eq:Main3-}) we conclude that the preferred end is to be found among the ones with bigger mass, because as it is easy to see an end with more positive mass do better than an end of less mass when we are considering large volumes. So from this perspective the worst case is the one considered in Theorem \ref{Res:Existence} in which all the masses $m_i$ are equal to their common value $m$ and in which we cannot say a priori which is the end that the isoperimetric regions for large volumes will prefer. However, Theorem \ref{Res:Existence} says that also in case of equal masses the number of ends in which the isoperimetric regions for large volumes concentrates is exactly one, but this end could vary from an isoperimetric region to another. An example of this behavior is given by Corollary $16$ of \cite{BrendleEichmair}, in which there are two ends and exactly two isoperimetric regions for the same large volume and they are obtained one from each other by reflection across the horizon, and each one of these isoperimetric regions chooses to have the biggest amount of mass in one end or in the other. 
\end{Rem}  
After this informal presentation of the proof of Corollary \ref{Cor:multimass}, we are ready to go into its details.

\begin{Dem} Here we treat the case in which the masses are not all equals, the case of equal masse being already treated in Theorem \ref{Res:Existence}. Without loss of generality we can assume that $1\in\mathcal{I}$, i.e., $$m_1=\max\{m_1,\dots,m_{\tilde{N}}\}.$$ We will prove the corollary by contradiction. To this aim, suppose that the conclusion of Corollary  \ref{Cor:multimass} is false, then there exists a sequence of isoperimetric regions $\Omega_j$ such that $V(\Omega_j)=v_j\rightarrow +\infty$, and $$E_{\Omega_j}\notin\{E_i\}_{i\in\mathcal{I}}.$$ Now we construct a competitor $\Omega'_j:=(\Omega_j\setminus E_{\Omega_j})\mathring{\cup}\tilde{B}^1_{r_j}$, such that $V(\tilde{B}^1_{r_j})=v'_j+v''_j$, with $v''_j:=V(\Omega_j\cap E_1)$ and $v'_j:=V(\Omega_j\setminus E_{\Omega_j})$. Roughly speaking it is like subtract the volume of $\Omega_j$ inside $E_{\Omega_j}$ and to put it inside the end $E_1$. As in the proof of Lemma \ref{Lemma:IsopregionsLargeVolume}, also in case of different masses we have that $v''_j$ is uniformly bounded and $v'_j\rightarrow+\infty$. By construction $V(\Omega'_j)=V(\Omega_j)=v_j$. Furthermore, it is not too hard to prove that we have the following estimates
\begin{equation}\label{Eq:Multimass}
A(\partial\Omega'_j)-A(\partial\Omega_j)\leq-\left(m^*_1-m^*_{E_{\Omega_j}}\right){v'_j}^{\frac{1}{n}}+o({v'_j}^{\frac{1}{n}}).
\end{equation}
This last estimate follows from an application of an analog of Lemma \ref{Lemma:IsopregionsLargeVolume} in case of different masses which goes mutatis mutandis and uses in a crucial way Theorem $4.1$ of \cite{EichmairMetzgerInv}. This cannot be avoided because again we need to control what happens to the area $A(\partial\Omega_j\cap E_{\Omega_j})$. The right hand side of (\ref{Eq:Multimass}), becomes strictly negative for $j\rightarrow+\infty$, since we have assumed $m^*_1-m^*_{E_{\Omega_j}}>0$. This yields to the desired contradiction. 
\end{Dem}  

Here we prove Theorem \ref{Res:Existence1}.   
   
      \begin{Dem} By Proposition 12 of \cite{BrendleEichmair} and equation (\ref{Eq:AsymptoticSchwarzschild2}) we get by a direct calculation that for a given $0<v<V(M)$ and any compact set $K\subseteq M$ there exists a smooth region $D\subset M\setminus K$ such that $V(D)=v$ and 
      \begin{equation}\label{Eq:Existence11}
                 A(\partial D)<c_nv^{\frac{n-1}{n}}=I_{\R^n}(v).
      \end{equation} 
      $D$ is obtained by perturbing the closed balls $B:=\{x: |x-a|\leq r\}$, for bounded radius $r$ and big $|a|$. The remaining part of the proof follows exactly the same scheme of Theorem $13$ of \cite{BrendleEichmair}, that was previously employed in another context in the proof of Theorem $1.1$ of \cite{MonNar}.  Now, using Theorem 1 of \cite{NarAsian}, reported here in Theorem \ref{Cor:Genexistence} we get that there exists a generalized isoperimetric region $\Omega=\Omega_1\mathring{\cup}\Omega_{\infty}$, both $\Omega_1\subseteq M$ and $\Omega_{\infty}\subseteq\R^n$ are isoperimetric regions in their own volumes in their respective ambient manifolds, with $V(\Omega)=v$, $V(\Omega_1)=v_1$, $V(\Omega_{\infty})=v_{\infty}$, $v=v_1+v_{\infty}$, moreover by Theorem 3 of \cite{NarAsian} $\Omega_1$ is bounded. If $\Omega_{\infty}=\emptyset$, the theorem follows promptly. Suppose, now that $\Omega_{\infty}\neq\emptyset$, one can chose as before a domain $D\subseteq M\setminus\Omega_1$ such that $V(D)=v_{\infty}$, $A(\partial D)<c_n{v_{\infty}}^{\frac{n-1}{n}}=I_{\R^n}(v_{\infty})$. This yields to the construction of a competitor $\Omega':=\Omega_1\mathring{\cup}D\subseteq M$ such that $V(\Omega')=v$ and $A(\partial\Omega')=A(\partial\Omega_1)+A(\partial D)<I_M(v)=A(\partial\Omega)$, this leads to a contradiction, hence $D_{\infty}=\emptyset$ and the theorem follows.       
      \end{Dem}
      \begin{Rem} As a final remark we observe that the hypothesis of convergence of the metric tensor stated in (\ref{Eq:AsymptoticSchwarzschild2}) are necessary for the proof of Theorem \ref{Res:Existence1}, because a weaker rate of convergence could destroy the estimate (\ref{Eq:Existence11}), when passing from the model Scharzschild metric to a $C^0$-asymptotically one.
      \end{Rem}  
      \newpage
      \markboth{References}{References}
      \bibliographystyle{alpha}
      \bibliography{these}
      \addcontentsline{toc}{section}{\numberline{}References}
      \emph{Stefano Nardulli\\ Departamento de Matem\'atica\\ Instituto de Matem\'atica\\ UFRJ-Universidade Federal do Rio de Janeiro, Brazil\\ email: nardulli@im.ufrj.br\\ \\ Abraham Henrique Mu\~noz Flores\\ Ph.D student\\Instituto de Matem\'atica\\UFRJ-Universidade Federal do Rio de Janeiro, Brazil\\email: abrahamemf@gmail.com} 
\end{document}